\documentclass[%
 reprint,
 amsmath,amssymb,
 aps,
]{revtex4-2}

\usepackage{graphicx}
\usepackage{dcolumn}
\usepackage{bm}
\usepackage{xcolor}

\begin{document}

\preprint{APS/123-QED}

\title{Ensemble-Based Global Search Framework for the Design Optimization of Fabrication-Constrained Freeform Devices}

\author{Seokhwan Min}
\author{Junhyung Park}%
\author{Jonghwa Shin}
 \email{qubit@kaist.ac.kr}
\affiliation{%
 Department of Materials Science and Engineering,\\
 Korea Advanced Institute of Science and Technology,\\
 Daejeon 34141, Republic of Korea
}%

\date{\today}

\begin{abstract}
Although freeform devices with complex internal structures promise drastic increases in performance, the discreteness of the set of available materials presents challenges for gradient-based optimization necessary for the efficient exploration of the high-dimensional freeform parameter space. Several schemes have been devised to utilize a continuous latent parameter space that maps to actual discrete designs, but none thus far simultaneously achieves full differentiability and strictly feasible material choices during optimization. Here, we propose the Gaussian ensemble gradient descent framework, which transforms the piecewise-constant fabrication-constrained cost function by convolution with a Gaussian kernel to render it differentiable. The transformed cost and gradient are estimated through ensemble sampling, which is combined with variance reduction methods that greatly improve the sampling efficiency in high-dimensional parameter spaces. Furthermore, the use of ensemble sampling within a gradient descent framework leads to the effective hybridization of the exploration and exploitation strengths of population- and gradient-based methods, respectively.
\end{abstract}

\maketitle


\section{\label{sec:introduction}Introduction}

The past several decades have seen groundbreaking improvements in multiple fields of engineering through freeform design for a plethora of applications including integrated photonic components~\cite{piggot15,su17,vercruysse19,lalau-keraly13}, optical metasurfaces~\cite{sell17_1,sell17_2,jiang19,christiansen19,zeng21,zhao21,lee24,kim24_1,miller12,liang13}, heat sinks and exchangers~\cite{li16,cohen22,fawaz22,erfani24,ali24,romano22}, acoustic metasurfaces~\cite{xu20,emoto23,fu25,cheng23,kurioka23}, and mechanical devices and metamaterials~\cite{aage17,wang25,li24,deaton16,zhu16,zhu20}. Such improvements come from the greatly increased degrees of freedom of pixelated freeform devices that lead to highly non-intuitive designs unlike simpler geometries such as cylinders or cuboids.

\begin{figure*}
\includegraphics[width=0.7\textwidth]{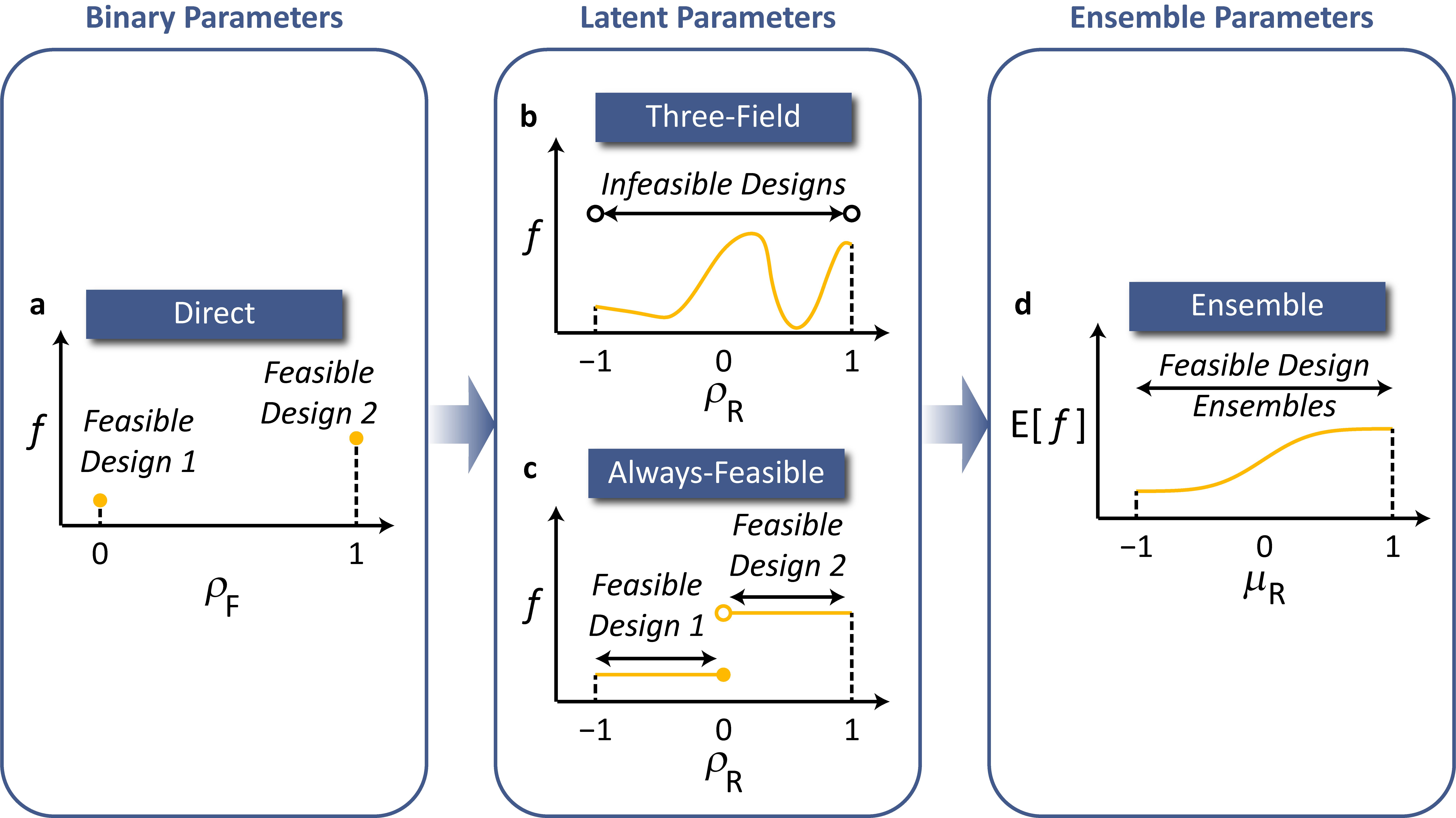}
\caption{\label{fig:principle} A summary of freeform design parameterization schemes used under previous and currently proposed optimization frameworks. \textbf{a}, A direct parameterization of densities corresponding to feasible designs lead to a cost function with discrete domain and range, which presents difficulties for gradient descent. \textbf{b},\textbf{c}, Latent parameters are used in more advanced frameworks such as the three-field method and the always-feasible parameterization to make the domain continuous. However, feasible designs only constitute a small fraction of the parameter space under the three-field parameterization (\textbf{b}), while the cost function is not differentiable everywhere under the always-feasible parameterization (\textbf{c}). \textbf{d}, Under the ensemble parameterization, the always-feasible cost function is probed using a Gaussian PDF whose mean is the optimization parameter. Every value for the mean $\mu_\text{R}$ represents a unique ensemble of feasible designs, and the objective is to minimize the average cost. Thus, both the domain and range of the optimization function can be made continuous without compromising design feasibility.}
\end{figure*}

Nevertheless, the effective design of fabrication-friendly freeform designs remains a challenge. Many common fabrication methods (e.g., photolithography) limit the constituents of the fabricated devices to two different materials (one of which is often empty space or an encapsulating material). In addition, fabrication equipments have associated resolutions that define the minimum feature sizes of the fabricated device. For pixel-parameterized designs, the former represents a pixel-wise constraint that restricts the parameter value at each pixel to values of 0 (material A) or 1 (material B). The latter is an inter-pixel constraint that enforces a minimum size to clusters of neighboring pixels of the same material. We will refer to designs that conform to both constraints as feasible designs. The discrete nature of both the domain and the range of the cost function caused by the material constraint precludes a direct exploration of the feasible design space via gradient descent (Fig.~\ref{fig:principle}a), while the high dimensionality of the freeform parameter space greatly reduces the efficiency of gradient-free heuristic approaches such as particle swarm optimization (PSO) and genetic algorithm. Thus, specialized parameterization techniques are necessary to facilitate effective optimization of feasible designs.

Initial approaches used continuous densities that interpolated between the relevant properties of materials A and B to enable gradient-based optimization, but such grayscale designs have limited fabrication feasibility. Subsequent developments led to a multi-level parameterization where continuous latent parameters are used to indirectly represent the actual design (Fig.~\ref{fig:principle}b--c). These approaches can be broadly classified into two categories: differentiable parameterizations that map latent parameters to near-feasible designs and non-differentiable parameterizations that map latent parameters to fully feasible designs. The former category is exemplified by the three-field method in which the latent density is filtered (blurred) using a kernel of choice and projected towards 0 or 1 through a hyperbolic tangent function whose projection strength is increased every $S$ iterations~\cite{wang11,lazarov16}. Some variations also add regularization~\cite{sell17_1,sell17_2} or inequality constraints that activate near the end~\cite{hammond21}. Although these methods are differentiable, they do not prevent optimization from progressing towards local optima in infeasible regions (Fig.~\ref{fig:principle}b). Therefore, increases in the design feasibility during and at the end of optimization tend to act against the cost function gradient, resulting in trade-offs in device performance~\cite{sell17_1,christiansen19,kim24_1,jin20}.

The latter category of latent parameterizations was pioneered by Schubert \textit{et al}. who proposed the always-feasible parameterization scheme that maps latent parameters to fully feasible designs at every iteration using a feasible design generator (FDG)~\cite{schubert22}. The cost is always computed with respect to the feasible design, preventing convergence to infeasible designs during optimization (Fig.~\ref{fig:principle}c). Despite this advantage, the FDG is non-differentiable and low-accuracy straight-through estimators (STE) are used as substitute gradients~\cite{schubert22,schubert25}, which degrades the quality of the optimized designs. Recently, Chen \textit{et al}. have proposed a differentiable near-feasible design generator that transforms latent designs into almost-feasible designs which can then be thresholded to fully feasible designs~\cite{chen24}. The close resemblance between the almost-feasible and fully feasible designs allows the gradient of the thresholding step to be better approximated by STE. Nevertheless, strictly feasible designs still remain non-differentiable even under this framework. In summary, latent parameterization requires one to choose between differentiability and feasibility, and to our knowledge, no parameterization scheme has been reported thus far that simultaneously allows for both during optimization.

In this work, we propose an alternative Gaussian ensemble gradient descent (GEGD) framework that combines the differentiability of the three-field method and the feasibility of the always-feasible parameterization. Our framework maintains an ensemble of feasible designs sampled by a multivariate Gaussian distribution in the latent density parameter space. The optimization objective is to find the mean of the sampling distribution that minimizes the expected cost (Fig.~\ref{fig:principle}d). The expected cost is equivalent to the convolution of the original cost function with a Gaussian probability density function (PDF). Therefore, both its domain and range are continuous, and it is always differentiable with respect to the distribution mean as long as the original cost function is finite with well-defined integrals. In actual implementation, the true expected cost is estimated by an ensemble average of Monte Carlo samples. The concept of minimizing the expected cost has previously been used in several general-purpose heuristic algorithms such as natural evolution strategies~\cite{wierstra14}, but they have seen limited use in freeform topology optimization. On the other hand, GEGD is specifically designed for effectively navigating the non-differentiable, high-dimensional  parameter space of freeform designs. Importantly, GEGD uses momentum-based parameter updates, radial basis function (RBF) sampling covariance, and control variates~\cite{mohamed20,jakeman23} that drastically improve the sampling efficiency for the estimation of the expected cost and its gradient. Through this, we make population-based optimization computationally tractable even in the high-dimensional parameter spaces typical of freeform devices. Furthermore, we show that the non-local cost function information gathered by ensemble sampling allows the effective hybridization of population-based global search and gradient-based convergence. Through benchmark tests on two representative nanophotonic design problems, we demonstrate the superior performance of our framework compared to conventional gradient and population-based methods. We note that despite our focus on nanophotonic benchmark problems, our framework is applicable for the design of density-based freeform designs in any field of engineering.

\section{\label{sec:results}Results}

\subsection{\label{sec:ensemble_parameterization}Ensemble Parameterization}

\begin{figure*}
\includegraphics[width=0.7\textwidth]{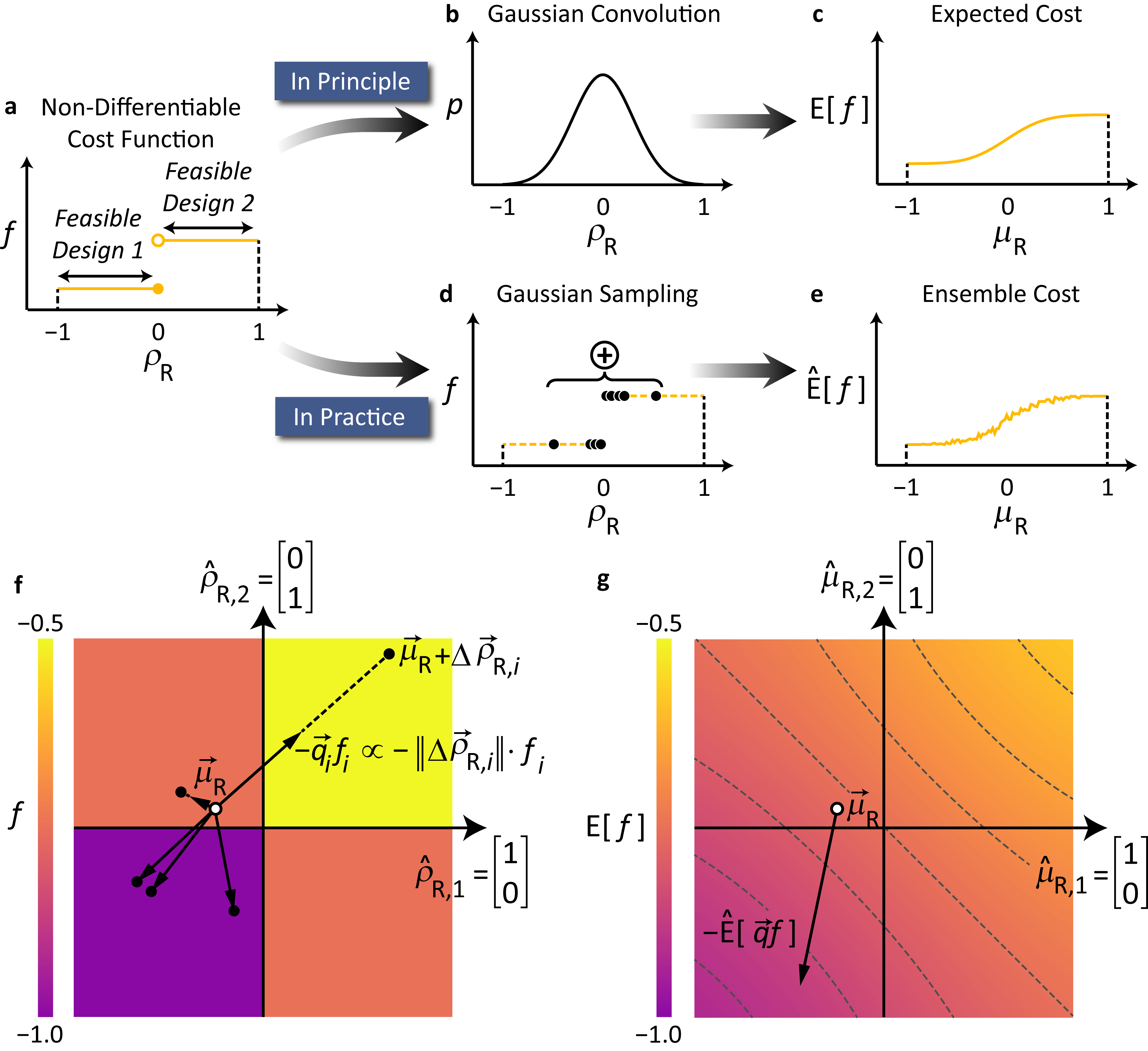}
\caption{\label{fig:implementation}
Implementation details regarding the ensemble cost and gradient evaluations. \textbf{a--e}, Given a non-differentiable cost function (\textbf{a}), convolution with a Gaussian (\textbf{b}) yields a reparametrized cost function (equivalent to the expectation of the original cost sampled under the Gaussian distribution) that is differentiable with respect to the Gaussian distribution parameters (\textbf{c}). In practice, Monte Carlo sampling (\textbf{d}) is done on the original cost to estimate the expected cost (\textbf{e}). \textbf{f}, The ensemble gradient is computed by averaging unit vectors pointing from the mean towards sampled points, each scaled by the respective sampled cost and the distance of the sampled point from the mean. \textbf{g}, Contour plot of the expected cost showing the ensemble gradient computed from the samples in \textbf{f}.}
\end{figure*}

We first briefly describe the FDG proposed by Schubert \textit{et al}. to facilitate the subsequent discussion of our framework. A latent density distribution $\vec{\rho}_\text{L}$, $\left( \rho_{\text{L},i} \in \left[-1,1\right] \right)$ is filtered and projected into a reward matrix $\vec{\rho}_\text{R}$, $\left( \rho_{\text{R},i} \in \left[-1,1\right] \right)$ that is provided as input to the FDG. The FDG makes solid and void `touches' with a predefined circular brush on an empty design. Their sequence is based on the `score' values computed by convolving the reward matrix with the brush centered on the location of each touch. The size of the brush determines the minimum feature size and the minimum radius of curvature of the resulting feasible design $\vec{\rho}_\text{F}$, $\left( \rho_{\text{F},i} \in \left\{0,1\right\} \right)$. If at any point an undetermined pixel becomes impossible to reach using a solid (void) brush without flipping void (solid) pixels, its density is automatically assigned to 0 (1). This process results in a piecewise constant cost function $f\left(\vec{\rho}_\text{R}\right)$ in which all $\vec{\rho}_\text{R}$ within a given region map to the same $\vec{\rho}_\text{F}$, with abrupt discontinuities at the boundaries (Fig.~\ref{fig:implementation}a).

Our GEGD framework replaces $\vec{\rho}_\text{L}$ and $\vec{\rho}_\text{R}$ with the mean latent density $\vec{\mu}_\text{L}$ and the mean reward matrix $\vec{\mu}_\text{R}$, respectively. The latter defines the center of a multivariate Gaussian PDF, $p$ (Fig.~\ref{fig:implementation}b), which is convolved with the optimization cost function as follows:
\begin{equation}
\label{eq:gaussian_convolution}
    f' \left( \vec{\mu}_\text{R} \right) = \\
    \int p \left( \vec{\rho}_\text{R}, \vec{\mu}_\text{R}, \sigma_\text{R}^2 \boldsymbol{\Sigma} \right) \\
    \cdot f \left( \vec{\rho}_\text{R} \right) \\
    \text{d} \vec{\rho}_\text{R}
\end{equation}
where $\vec{\mu}_\text{R}$ and $\sigma_\text{R}^2 \boldsymbol{\Sigma}$ are the mean and covariance of the multivariate Gaussian distribution. As long as $f$ is finite with well-defined integrals, the function $f'$ is differentiable with respect to $\vec{\mu}_\text{R}$ (Fig.~\ref{fig:implementation}c). Therefore, instead of solving for $\vec{\rho}_\text{R}$ that minimizes $f$, the optimization problem can be reformulated as solving for $\vec{\mu}_\text{R}$ that minimizes $f'$.

Equation (\ref{eq:gaussian_convolution}) is challenging to evaluate exactly for typical electromagnetic cost functions, which often lack closed-form expressions in terms of the permittivity distribution. Instead, it can be estimated using the fact that it is conceptually equivalent to the expectation of $f$ under Gaussian sampling. Then, the true expectation can be estimated from a finite number of $M$ samples taken from $f$ (Fig.~\ref{fig:implementation}d--e):
\begin{subequations}
    \label{eq:ensemble_cost}
    \begin{equation}
        f' \left( \vec{\mu}_\text{R} \right) \\
        \approx \left< f \left( \vec{\mu}_\text{R} + \Delta \vec{\rho}_{\text{R},i} \right) \right>
    \end{equation}
    \begin{equation}
        \Delta \vec{\rho}_{\text{R},i} \\
        \sim N \left( \vec{0}, \sigma_\text{R}^2 \boldsymbol{\Sigma} \right)
    \end{equation}
\end{subequations}
where $\left< ... \right>$ indicates an ensemble average. We will refer to this estimate as the “ensemble cost”.

For optimization, we need to compute the expected cost gradient with respect to $\vec{\mu}_\text{R}$. Differentiating Eq.~(\ref{eq:gaussian_convolution}) yields:
\begin{subequations}
\label{eq:ensemble_gradient}
    \begin{equation}
    \label{subeq:expected_cost_gradient}
        \frac{\text{d} f'}{\text{d} \vec{\mu}_\text{R}} \\
        = \int p \left( \vec{\rho}_\text{R}; \vec{\mu}_\text{R}, \sigma_\text{R}^2 \boldsymbol{\Sigma} \right) \\
        \cdot \vec{q} \left( \vec{\rho}_\text{R}, \vec{\mu}_\text{R} \right)\\
        \cdot f \left(\vec{\rho}_\text{R} \right) \text{d} \vec{\rho}_\text{R} \\
    \end{equation}
    \begin{equation}
    \label{subeq:ensemble_gradient}
        \frac{\text{d} f'}{\text{d} \vec{\mu}_\text{R}} \\
        \approx \left< \vec{q} \left( \vec{\mu}_\text{R} + \Delta \vec{\rho}_{\text{R},i}, \vec{\mu}_\text{R} \right) \cdot f \left( \vec{\mu}_\text{R} + \Delta \vec{\rho}_{\text{R},i} \right) \right>
    \end{equation}
    \begin{equation}
        \vec{q} \left( \vec{\rho}_\text{R}, \vec{\mu}_\text{R} \right) = \sigma_\text{R}^{-2} \boldsymbol{\Sigma}^{-1} \\
        \left( \vec{\rho}_\text{R} - \vec{\mu}_\text{R} \right)
    \end{equation}
\end{subequations}
The integral in Eq.~(\ref{subeq:expected_cost_gradient}) is equivalent to the expectation of a new vector-valued function $\vec{q} f$ under the same Gaussian distribution used to evaluate the expected cost. The expectation, therefore, can also be estimated using an ensemble average in the same manner as the expected cost itself. We will refer to this estimate as the “ensemble gradient” (Eq.~(\ref{subeq:ensemble_gradient})). Crucially, Eq.~(\ref{eq:ensemble_gradient}) does not involve the gradient of the original electromagnetic cost function $f$, which means it does not require adjoint simulations or autogradable simulators. In fact, $f$ need not be differentiable at all. In addition, the computation of the ensemble gradient does not incur any additional simulations, because it can be computed using the same sampled costs used for the computation of the ensemble cost (Eq.~(\ref{eq:ensemble_cost})). The ensemble gradient is then backpropagated through the filtering and projection operations with analytically known gradients to obtain the overall ensemble gradient with respect to $\vec{\mu}_\text{L}$.

\begin{figure*}
\includegraphics[width=0.8\textwidth]{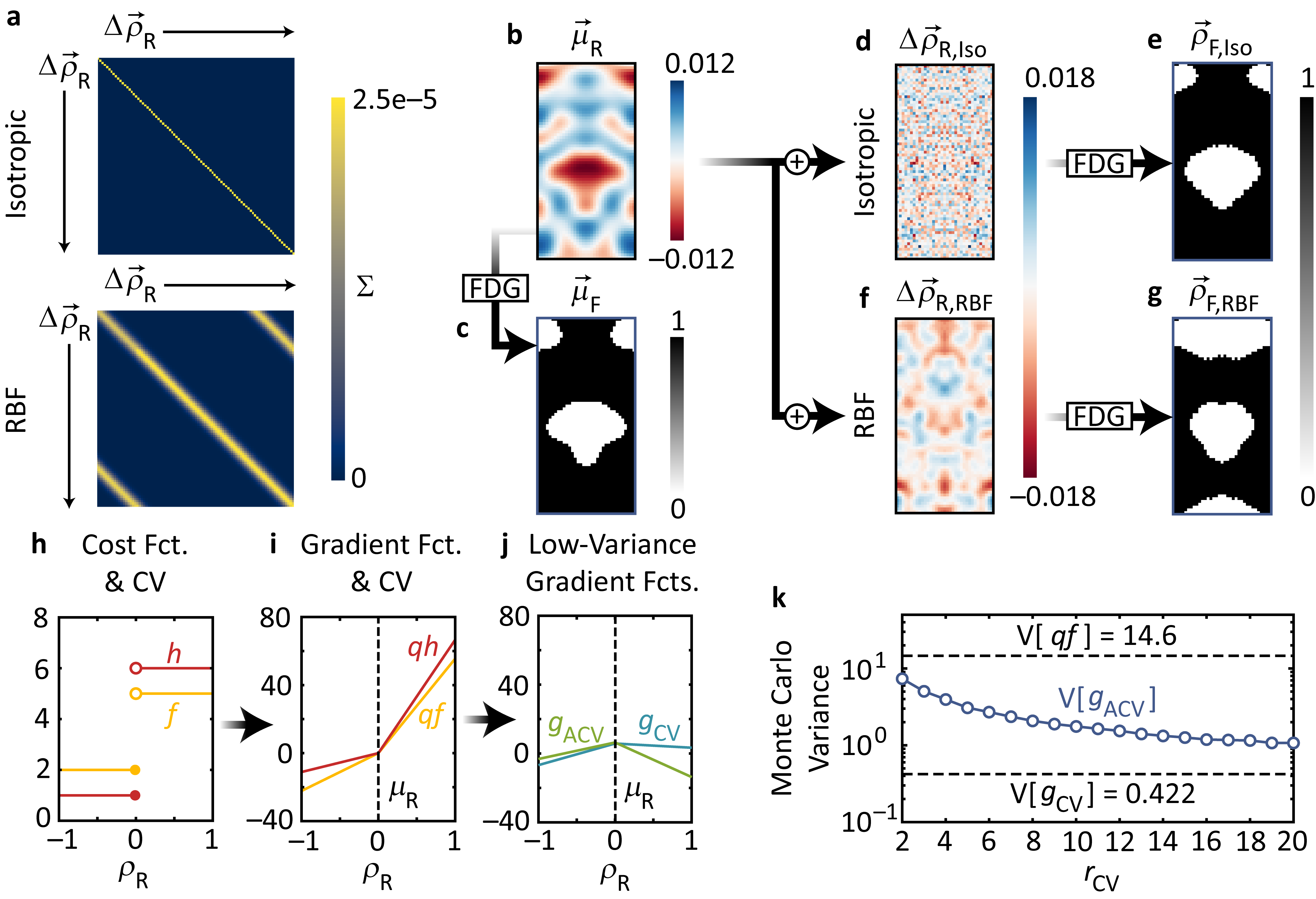}
\caption{\label{fig:sampling_efficiency} Methods for improving Monte Carlo sampling efficiency for freeform designs with large numbers of parameters. \textbf{a}, Examples of isotropic (top) and RBF (bottom) covariance matrices. The plots shown are $88 \times 88$ submatrices of the full $1260 \times 1260$ covariance matrices that produce the perturbations in \textbf{d} and \textbf{f}. \textbf{b},\textbf{c}, An example mean reward matrix (\textbf{b}) and the corresponding feasible design (\textbf{c}). \textbf{d},\textbf{e}, Isotropic perturbation (\textbf{d}) and subsequent feasible design generation (\textbf{e}) leads to a design that is only slightly different from the original design. \textbf{f},\textbf{g}, RBF perturbation (\textbf{f}) results in significant topological changes in the feasible design (\textbf{g}). \textbf{h--j}, Plots describing the construction of low-variance gradient functions which are sampled to compute the ensemble gradient. Original ($f$) and control variate ($h$) cost functions (\textbf{h}) are multiplied by $g_\mu$ to construct the gradient function ($g_\mu f$) and its CV ($g_\mu h$) (\textbf{i}). Low-variance gradient functions are constructed by subtracting a scaled and shifted $g_\mu h$ from $g_\mu f$ (\textbf{j}). The approximate CV gradient function was computed using 10 high-fidelity and 100 low-fidelity samples ($r_\text{CV} = 10$). \textbf{k}, Monte Carlo sampling variance as a function of the low-to-high-fidelity sampling ratio ($r_\text{CV}$). The upper and lower limits represents variance for the original gradient function and the exact CV gradient function, respectively.}
\end{figure*}

Intuitively, the ensemble average in Eq.~(\ref{subeq:ensemble_gradient}) can be interpreted as a weighted average of directions $\Delta \vec{\rho}_\text{R} / \left\Vert \Delta \vec{\rho}_\text{R} \right\Vert$ sampled around $\vec{\mu}_\text{R}$. Assuming a scalar covariance for simplicity (i.e., $\boldsymbol{\Sigma} = \textbf{I}$), the weights are proportional to the sampled cost multiplied by the distance of the sampled point from the distribution center.
\begin{equation}
\label{eq:ensemble_gradient_interpretation}
    \frac{\text{d} f'}{\text{d} \vec{\mu}_\text{R}} \\
    \approx \left< \frac{1}{\sigma_\text{R}^2} \\
    f \left( \vec{\rho}_\text{R} \right) \\
    \left\Vert \Delta \vec{\rho}_\text{R} \right\Vert \\
    \cdot \frac{\Delta \vec{\rho}_\text{R}}{\left\Vert \Delta \vec{\rho}_\text{R} \right\Vert} \right>
\end{equation}
Figure~\ref{fig:implementation}f shows a simplified example in two dimensions in which the ensemble gradient is computed using five samples. As can be seen in Eq.~(\ref{eq:ensemble_gradient_interpretation}), more weight is given to samples that have large absolute costs. It also seemingly indicates that samples far away from the mean will dominate the gradient, but the sampling probability itself falls exponentially away from the mean. The product of the sample distance and the sampling probability is $\propto \left\Vert \Delta \vec{\rho}_\text{R} \right\Vert \text{exp} \left( -\left\Vert \Delta \vec{\rho}_\text{R} \right\Vert^2 / 2 \sigma_\text{R}^2 \right)$ whose extrema lie at $\left\Vert \Delta \vec{\rho}_\text{R} \right\Vert = \pm \sigma_\text{R}$. Therefore, samples that are approximately one standard deviation away from $\vec{\mu}_\text{R}$ tend to have the largest representation in the ensemble gradient. This means that the optimizer can in principle hop between local minima spaced apart by approximately one standard deviation or less.

\subsection{\label{sec:sampling_efficiency}Monte Carlo Sampling Efficiency Enhancement Techniques}

Monte Carlo sampling is notoriously noisy in high-dimensional parameter spaces (on the order of $10^2$ – $10^4$) typical for freeform devices. In addition, electromagnetic simulations needed for each sampled cost tend to be computationally expensive. Therefore, we employ three different strategies to improve sampling efficiency and reduce ensemble gradient estimation noise.

Firstly, we used momentum-based ADAM updates of the mean latent density, which stabilize noisy gradients through exponential moving averages. The implementation details are described in ``Methods".

Secondly, the cost function was sampled using an anisotropic Gaussian distribution with a radial basis function (RBF) covariance that has distance-based correlation between pixels (Fig.~\ref{fig:sampling_efficiency}a).
\begin{equation}
\label{eq:RBF_covariance}
    \boldsymbol{\Sigma}_{\text{RBF},ij} \\
    = \text{exp} \left( -\left\Vert \vec{x}_i - \vec{x}_j \right\Vert^2 / \sigma_\text{RBF}^2 \right) \\
    + \delta_{ij} \epsilon
\end{equation}
where $\vec{x}_i$ and $\vec{x}_j$ are the coordinates of pixels $i$ and $j$, respectively, and $\sigma_\text{RBF} = L_\text{min} \sqrt{2} / 4$. $\delta_{ij}$ is the Kronecker delta and $\epsilon$ is a regularization factor that stabilizes the numerical inversion of $\boldsymbol{\Sigma}_\text{RBF}$ (see ``Methods"). The overall covariance is given by $\sigma_\text{R}^2 \boldsymbol{\Sigma}_\text{RBF}$. The perturbations $\Delta \vec{\rho}_\text{R}$ generated using the RBF covariance vary spatially on the order of $\sigma_\text{RBF} \propto L_\text{min}$. This leads to significant variations in the sequence of touches made by the FDG and the resulting $\vec{\rho}_\text{F}$ (Fig.~\ref{fig:sampling_efficiency}f--g), as opposed to perturbations from an isotropic covariance $\boldsymbol{\Sigma} = \textbf{I}$ which would waste samples exploring near-identical $\vec{\rho}_\text{F}$ (Fig.~\ref{fig:sampling_efficiency}d--e).

Finally, we used approximate control variates to reduce Monte Carlo sampling variance for the ensemble gradient. The central ideal is to replace the original vector function $\vec{q} f$ with a substitute function $\vec{g}_\text{CV}$ with lower variance~\cite{mohamed20,jakeman23} (Fig.~\ref{fig:sampling_efficiency}h--j):
\begin{subequations}
\label{eq:exact_CV}
\begin{equation}
\label{subeq:gcv}
    \vec{g}_\text{CV} \left( \vec{\rho}_\text{R} \right) = \\
    \vec{q} f \\
    -\beta_\text{CV} \left\{ \vec{q} h - \text{E} \left[ \vec{q} h \right] \right\}
\end{equation}
\begin{equation}
\label{subeq:beta_CV}
    \beta_\text{CV} \\
    = \overline{\text{Cov}} \left[ \vec{q} f, \vec{q} h \right] \\
    / \overline{\text{V}} \left[ \vec{q} h \right]
\end{equation}
\begin{equation}
\label{subeq:variance_reduction_exact}
    \frac{\overline{\text{V}} \left[ \vec{g}_\text{CV} \right]}{\overline{\text{V}} \left[ \vec{q} f \right]} \\
    = 1 - \overline{\text{Corr}} \left[ \vec{q} f, \vec{q} h \right]^2
\end{equation}
\end{subequations}
$\overline{\text{Q}} \left[ \vec{v} \right]$ denotes the average of the quantity Q over all vector components of $\vec{v}$. The function $h$, the control variate (CV), is a function that closely approximates the original function $f$ (Fig.~\ref{fig:sampling_efficiency}h). If $h$ is analytically known such that $\text{E} \left[ \vec{q} h \right]$ can be evaluated exactly, the low-variance gradient function $\vec{g}_\text{CV}$ (Eq.~(\ref{subeq:gcv})) has optimal variance reduction (Eq.~(\ref{subeq:variance_reduction_exact})). However, simple analytic models that correlate well with simulations are not readily available for most electromagnetic systems. To ensure good correlation, $f$ and $h$ can be taken as cost functions evaluated using high- and low-fidelity simulations (e.g., RCWA with different Fourier order truncations), even though both are non-analytic. In such cases, the method of approximate control variates (ACV) can be used where $\text{E} \left[ \vec{q} h \right]$ is approximated by sampling~\cite{jakeman23}:
\begin{subequations}
\label{eq:approx_CV}
\begin{equation}
\label{subeq:expectation_ACV}
    \text{E} \left[ \vec{q} h \right] \\
    \approx \frac{1}{r_\text{CV} M} \\
    \sum_{m = 1}^{r_\text{CV} M} \vec{q} h
\end{equation}
\begin{equation}
\label{subeq:varianc_reduction_approx}
    \frac{\overline{\text{V}} \left[ \vec{g}_\text{ACV} \right]}{\overline{\text{V}} \left[ \vec{q} f \right]} \\
    = 1 - \frac{r_\text{CV} - 1}{r_\text{CV}} \overline{\text{Corr}} \left[ \vec{q} f, \vec{q} h \right]^2
\end{equation}
\end{subequations}
Thus, a total of $M$ samples are taken from $f$, and  $r_\text{CV} M$ samples (including the previous $M$ sampling points for $f$) are taken from $h$, where $r_\text{CV} = 2, 3, ...$ determines the ratio between the two. The variance reduction provided by ACV (Eq.~(\ref{subeq:varianc_reduction_approx})) is dependent on $r_\text{CV}$ (Fig.~\ref{fig:sampling_efficiency}d), and the optimal variance reduction (i.e., that of $\vec{g}_\text{CV}$) is recovered in the limit of $r_\text{CV} \rightarrow \infty$. In practice, $10 \leq r_\text{CV} \leq 20$ is often enough to produce a noticeable reduction in variance given a good correlation ($\geq 0.9$) between the original function and its CV.

\begin{figure}
\includegraphics[width=0.5\textwidth]{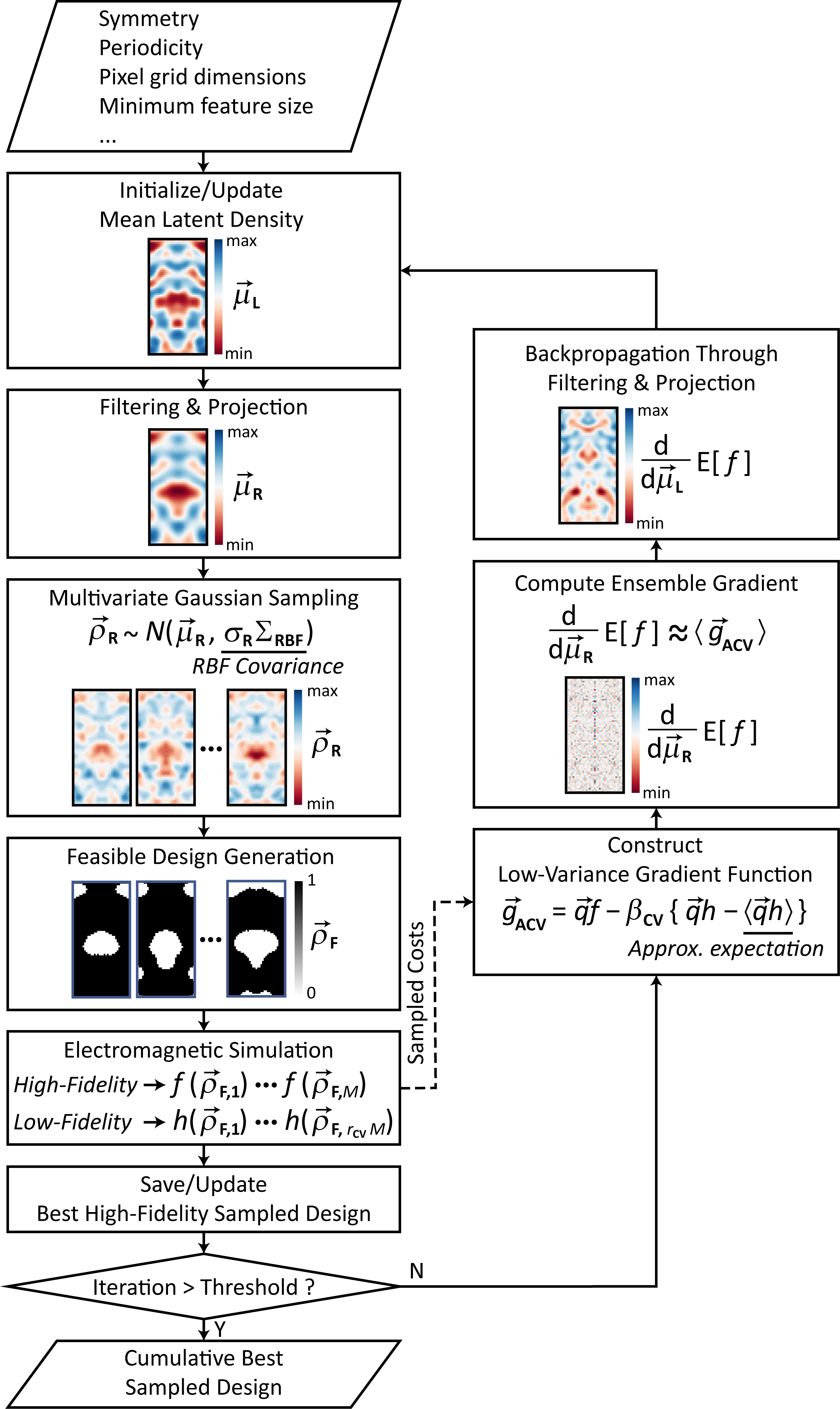}
\caption{\label{fig:flowchart} Flowchart of the overall Gaussian ensemble gradient descent framework. The mean latent density is filtered and projected to yield the mean reward matrix. This is used to generate a reward matrix ensemble, each of which is fed into a feasible design generator. The cost for each feasible design is simulated and ensemble averaged to yield an estimate for the expected cost. The gradient of the expected cost with respect to the mean reward matrix is estimated using the same sampled costs. This gradient is backpropagated through the filtering and projection operations to yield the final gradient. The best sampled cost and the corresponding feasible design is saved and updated in each iteration and provided as the final output after a set number of iterations.}
\end{figure}

Figure~\ref{fig:flowchart} shows the overall ensemble cost evaluation and gradient backpropagation scheme. The mean latent density $\vec{\mu}_\text{L}$ is the lowest-level parameter directly perturbed by the optimizer. Its filtered and projected form, the mean reward matrix $\vec{\mu}_\text{R}$, represents the mean of the Gaussian sampling distribution. A reward matrix ensemble is generated from the distribution and turned into feasible designs. The cost of each feasible design is evaluated, from which the ensemble cost is computed. The algorithm also stores and updates the best sampled design, which is what it outputs when the optimization terminates after a fixed number of iterations. Note that the optimized $\vec{\mu}_\text{L}$ is meaningful only in the context of its ability to increase the probability of sampling high-performance feasible designs. The ensemble gradient is computed by substituting the already sampled costs into the low-variance gradient function obtained through ACV.

\subsection{\label{sec:test_function}Radial Quasi-Invariance of Cost Functions Based on Feasible Design Generators}

In Sec.~\ref{sec:ensemble_parameterization}, we have described the capability of GEGD to hop between local minima using the non-local cost function information collected by ensemble sampling. However, in practice, $\sigma_\text{R}$ should be kept low so that the cost landscape of the ensemble does not become overly smoothed-out causing the ensemble gradient to be dominated by sampling noise. Thus, for reasonable $\sigma_\text{R}$, GEGD in itself is best described as a short-range global search algorithm that is somewhere between local methods such as gradient descent and long-range global methods such as PSO.

\begin{figure*}
\includegraphics[width=0.9\textwidth]{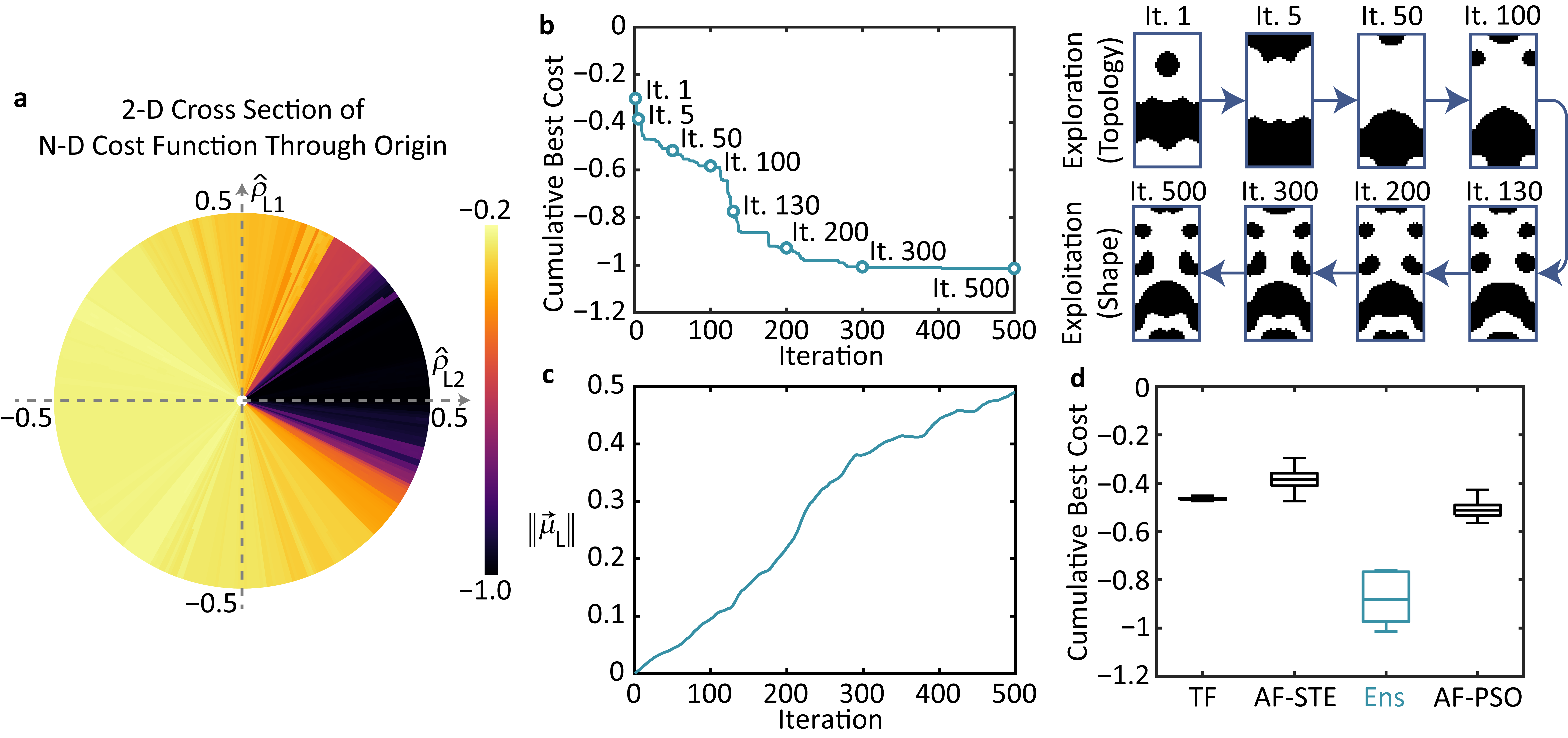}
\caption{\label{fig:test_function} Analytic test function benchmark results. \textbf{a}, A 2-dimensional cross section of the analytic cost function landscape through the origin along selected unit vectors $\hat{\rho}_{\text{L}1}$ and $\hat{\rho}_{\text{L}2}$. \textbf{b}, Cumulative best cost evolution during GEGD (left) with the corresponding feasible design at specific iterations (right). The algorithm performs topology exploration until approximately iteration 130, after which it fine tunes the shape of the obtained topology. \textbf{c}, The evolution of the Euclidean norm of the mean latent density during GEGD. \textbf{d}, Box-and-whisker plots of the best costs obtained using different optimization algorithms.}
\end{figure*}

However, a unique property of cost function landscapes involving the brush-based FDG allows even narrow sampling distributions to effectively cover a much larger set of feasible designs. As discussed in Sec.~\ref{sec:ensemble_parameterization}, the sequence of solid and void touches of the FDG depends on the score values, which may be scaled without changing the sequence. Scaling $\vec{\rho}_\text{L}$ roughly corresponds to scaling $\vec{\rho}_\text{R}$ and the score. Therefore, $\vec{\rho}_\text{F}$ is quasi-invariant to scaling $\vec{\rho}_\text{L}$ (i.e., ``radial quasi-invariance"), and the origin $\vec{\rho}_\text{L} = \vec{0}$ directly borders nearly the entire set of feasible designs in the parameter space (Fig.~\ref{fig:test_function}a). In such a cost function landscape, the best optimization strategy is to start at the origin and find the optimal direction to move along. Nevertheless, neither simple gradient descent nor heuristic algorithms can be used. The singular nature of the origin with undefined derivatives in all directions precludes gradient-based optimization, and the high-dimensional parameter space greatly reduces the efficiency of heuristic algorithms. However, GEGD can be initialized at the origin ($\vec{\mu}_\text{L}^{\left( 0 \right)}=\vec{0}$), because the Gaussian sampling smooths out the singularity.

To visualize the radial quasi-invariance, we define the following non-convex cost function on a design grid of $35 \times 70$ pixels.
\begin{equation}
    f_\text{test} = \\
    -\sum_{i=1}^{10} 3~\text{exp} \left[ -\frac{15}{N} \\
    \left\Vert \vec{\rho}_{\text{opt},i} - \vec{\rho}_\text{F} \left( \vec{\rho}_\text{L} \right) \right\Vert^2 \right]
\end{equation}
D1 symmetry (mirror symmetry across one axis) and a minimum feature size of 7 pixels were imposed on the device, leading to $N = 1260$ independent parameters. $\vec{\rho}_{\text{opt},i}$ are randomly generated filtered and projected grayscale designs that serve as local minima (Fig. S1). Figure~\ref{fig:test_function}a is a 2D cross section of the $N$-dimensional test function through the origin along selected unit vectors (Fig. S2), where the radial quasi-invariance can be clearly seen.

Figures~\ref{fig:test_function}b--c are results of a sample GEGD run performed on the test function. The initial optimization iterations ($< 130$) are spent near the origin exploring a wide range of different design topologies (Fig.~\ref{fig:test_function}b). Once a promising direction is identified, optimization accelerates away from the origin to fine-tune the shape of the design within a small angular slice of the parameter space. Due to the radial quasi-invariance of the cost function and the fact that $\sigma_\text{R}$ is kept constant throughout optimization, the ensemble gradient always has a component along the radial direction, and optimization always progresses outward from the origin (Fig.~\ref{fig:test_function}c).

We benchmarked the GEGD algorithm (including all variance reduction strategies) against three different conventional algorithms: the three-field parameterization coupled with the L-BFGS-B algorithm (TF), always-feasible optimization with straight-through estimator gradients (AF-STE), and particle swarm optimization with always-feasible parameterization (AF-PSO). Note that our implementation of grayscale optimization does not strictly enforce minimum feature size constraints. Refer to ``Methods" for algorithm implementation details and Supplementary Information S3 for a description of the CV used for GEGD.

For fair comparison, the number of different runs and the number of samples for the various algorithms were adjusted so that their overall computation costs were as close as possible. In terms of computation time, each iteration of GEGD and AF-PSO took the equivalent of 10 high-fidelity forward simulations. We assumed that the combined forward and backward simulations for TF and AF-STE take 1.5 times longer compared to forward-only simulations, so for every run of GEGD and AF-PSO, $10 / 1.5 \approx 7$ runs with different initial designs were performed for TF and AF-STE.

Figure~\ref{fig:test_function}d is a box-and-whisker plot showing the distributions of the costs obtained using all benchmark algorithms, where it is clear that GEGD outperforms all other algorithms by a significant margin. The fact that the algorithm was able to consistently identify high-performance designs is remarkable considering the large discrepancy between the number of samples (10 in terms of the effective computation cost) and the dimension of the parameter space (1260). This performance enhancement attests the effectiveness of our sample efficiency enhancement techniques (Fig. S3) and the hybridization of population-based global search and gradient-based local convergence.

\subsection{\label{sec:nanophotonic_benchmarks}Optimization Benchmarks for Nanophotonic Design Problems}

\begin{figure*}
\includegraphics[width=0.9\textwidth]{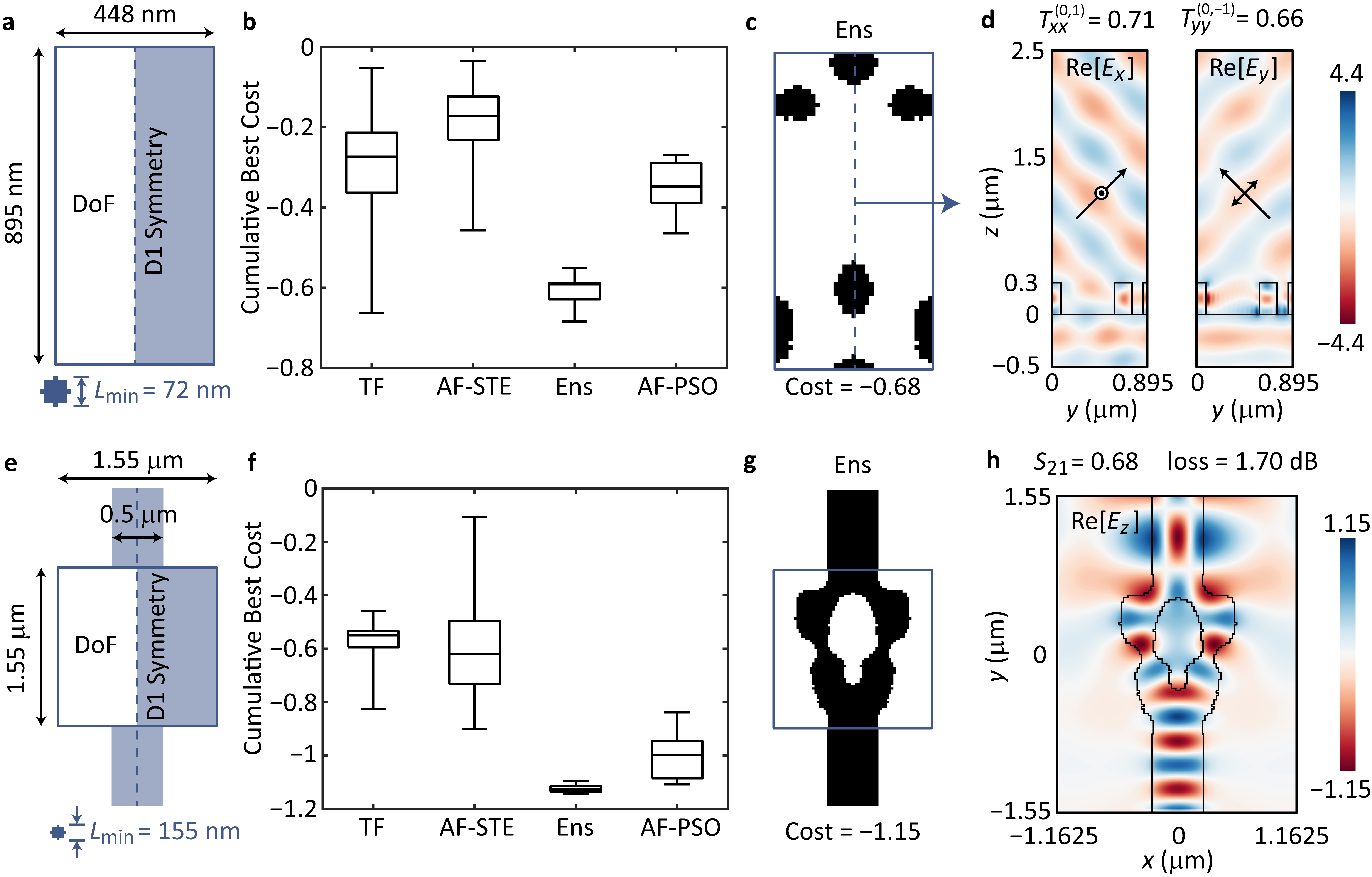}
\caption{\label{fig:nanophotonic_optimization} Nanophotonic optimization benchmarks. \textbf{a--d}, Polarization beamsplitter optimization. \textbf{a}, Geometric specifications. The unshaded region denotes the optimization degrees of freedom as specified by D1 symmetry. The minimum allowed feature size is shown at the bottom. \textbf{b}, Best costs obtained using different benchmark algorithms. \textbf{c}, Best design obtained by GEGD. \textbf{d}, $\text{Re} \left[ E_x \right]$ (left) and $\text{Re} \left[ E_y \right]$ (right) for $x$ and $y$-polarized incidence, respectively. The fields are plotted across the $yz$-plane at the center of the unit cell denoted by the dotted line in \textbf{c}. \textbf{e--h}, 2D mode converter optimization. \textbf{e}, Geometric specifications. The unshaded region denotes the optimization degrees of freedom as specified by D1 symmetry. The minimum allowed feature size is shown at the bottom. A 0.5 $\mu m$-wide waveguide is connected to either ends of the design region. \textbf{f}, Best costs obtained using different benchmark algorithms. \textbf{g}, Best design obtained by GEGD. \textbf{h}, $\text{Re} \left[ E_z \right]$ for $\text{TE}_0$ mode input from the bottom waveguide.}
\end{figure*}

In this section, we benchmark the four algorithms (TF, AF-STE, GEGD, and AF-PSO) on two different nanophotonic design problems: the design of polarization beamsplitters and 2D integrated mode converters.

The polarization beamsplitter was designed to diffract $x$- and $y$-polarized light toward $+45^\circ$ and $-45^\circ$, respectively, at a target wavelength of 633 nm. The design region has D1 symmetry and is divided into a $35 \times 70 \times 1$ array of pixels (Fig.~\ref{fig:nanophotonic_optimization}a). The optimization cost function was defined as:
\begin{equation}
    f \left( \vec{\rho}_\text{L} \right) \\
    = -\frac{1}{2} \left[ T_{xx}^{\left( +1 \right)} \left( \vec{\rho}_\text{L} \right) \\
    + T_{yy}^{\left( -1 \right)} \left( \vec{\rho}_\text{L} \right) \right]
\end{equation}
$T_{xx}^{\left( +1 \right)}$ and $T_{yy}^{\left( -1 \right)}$ denote the $x$-to-$x$ transmission towards $+45^\circ$ and $y$-to-$y$ transmission towards $-45^\circ$, respectively. Simulations were performed with TORCWA~\cite{kim23}. We used $10 \times 20$ harmonics for the high-fidelity simulations and $4 \times 8$ harmonics for the low-fidelity CV. The simulation times differ by a factor of $\approx 33$ using 32-thread multi-threading. The total computation cost per iteration for all algorithms was adjusted to be as close as possible to $10 t_\text{HF}$, where $t_\text{HF}$ is the runtime for a single high-fidelity forward simulation. Autodifferentiation-enabled simulations take $\approx 1.5$ times longer than forward-only simulations, so for every run of GEGD and AF-PSO, 7 runs were performed for TF and AF-STE.

The 2D mode converter was designed for operation at 1.55 $\mu$m where it converts the incoming $\text{TE}_0$ mode to $\text{TE}_2$. The design region has D1 symmetry and is divided into a $70 \times 70 \times 1$ array of pixels (Fig.~\ref{fig:nanophotonic_optimization}d). The cost was defined as follows:
\begin{equation}
    f \left( \vec{\rho}_\text{L} \right) \\
    = -\frac{\left| \int \vec{E}_\text{out} \left( \vec{\rho}_\text{L} \right) \cdot \vec{E}_\text{TE2}^*~\text{d}x \right|^2}{\left| \int \vec{E}_\text{TE2} \cdot \vec{E}_\text{TE2}^*~\text{d}x \right|^2}
\end{equation}
$\vec{E}_\text{out}$ and $\vec{E}_\text{TE2}$ denote the measured output field profile and the target $\text{TE}_2$ mode profile, respectively. The devices were simulated using Ansys FDTD where we used a mesh density of 5 per pixel side length for the high-fidelity simulations and 2 for the low-fidelity CV. The simulation times differ by a factor of $\approx 3.6$ using 12-thread multi-threading. The total computation cost per iteration for all algorithms was adjusted to be near $20 t_\text{HF}$. Adjoint simulations for TF and AF-STE take approximately the same time as the forward simulations, so for every run of GEGD and AF-PSO, 10 runs were performed for TF and AF-STE.

As shown in Figs.~\ref{fig:nanophotonic_optimization}b and f, GEGD outperforms all other algorithms by a significant margin for both benchmark problems. The best designs obtained using GEGD are shown in Figs.~\ref{fig:nanophotonic_optimization}c and g with the corresponding field plots in Figs.~\ref{fig:nanophotonic_optimization}d and h that confirm that the devices operate as specified. Additional data showing the evolution of the cost during optimization and the best designs for each algorithm can be found in Figs. S4--S5.

\section{\label{sec:discussion}Discussion}

We have proposed and benchmarked the GEGD algorithm on several nanophotonic design problems where it demonstrated superior performance compared to conventional methods even with the same computational cost. Our framework reformulates the optimization objective using ensemble parameters that define a multivariate Gaussian sampling distribution in the latent density space, whose convolution with the discontinuous landscape of the AF parameter space renders it smooth and differentiable. We also presented a method for estimating the gradient of the smoothed cost function using Monte Carlo sampling, whose sampling efficiency was enhanced through a combination of momentum updates, RBF covariance, and approximate CV. Finally, we demonstrated that the non-local cost function landscape information collected through Gaussian sampling imparts short-range minima hopping properties to GEGD. This is further enhanced by the radial quasi-invariance of the cost function landscape when using the brush-based FDG of Schubert \textit{et al}., which allows nearly all feasible designs to potentially be sampled near the origin even with narrow sampling distributions.

Our method is a general framework that may be applied to the design of a diverse range of density-based freeform devices, even outside of nanophotonics. Furthermore, it enables the adoption of efficient gradient-based optimization algorithms regardless of the differentiability of the cost function, because the ensemble parameterization bypasses the gradients of both the feasible design generator and the cost function. This makes GEGD potentially useful in various scenarios where freeform optimization was traditionally impossible or limited due to non-differentiable cost functions. For example, surface-enhanced Raman spectroscopy benefits from enhancing the volume-averaged near-field intensity in the void regions of a freeform design, which is challenging to directly optimize, since the void region itself is constantly changing during optimization. Further examples of non-differentiable cost functions may also be found in topological photonics, which involve non-differentiable topological invariants. We therefore expect GEGD to have a widespread impact in many applications in nanophotonics and beyond.

\section{\label{sec:methods}Methods}

\subsection{\label{sec:algorithm_details}Algorithm Implementation Details}

All four benchmark algorithms involve the filtering and projection of latent densities to either grayscale devices (TF) or reward matrices (AF-STE, AF-PSO, GEGD). Filtering was performed using a Gaussian kernel with standard deviation $\sigma_\text{f} = L_\text{min}$ (TF) or $\sqrt{2} L_\text{min} / 4$ (AF-STE, AF-PSO, GEGD). The resulting density values were projected using a hyperbolic tangent function with projection strength $\beta_\text{proj} = 8, 16, 32, 64, 128$ (TF) or 8 (AF-STE, AF-PSO, GEGD).

AF-STE, AF-PSO, and GEGD involve an additional feasible design generation step. For this, we used the brush-based algorithm by Schubert \textit{et al}.~\cite{schubert22}. The algorithm was adapted to include symmetry constraints by (1) ensuring that the reward matrix is symmetric and (2) making each touch simultaneously at multiple positions as dictated by the particular symmetry.

TF was coupled with the L-BFGS-B implementation in SciPy for parameter updates. The hyperbolic tangent projection strength was increased every 100 (test function, polarization beamsplitter) or 50 (mode converter) iterations, but slight deviations in the number of iterations were allowed as required by the line-search subroutine in L-BFGS-B.

AF-PSO used a variant of the PSO algorithm that probabilistically randomizes a subset of particle velocities every iteration (`Craziness')~\cite{dimoupoulos07}. The cognitive and social coefficients were set to 1.49 to balance exploration and exploitation. The inertia was initially set to 0.9 and was reduced by a factor of 0.95 whenever the global best did not update for 5 or more iterations. The craziness probability was set to 0.22 and 10\% of the particle velocities were reset if the probability was met.

ADAM parameter updates were used for AF-STE and GEGD. $\beta_1 = 0.667$, $\beta_2 = 0.9$, and $\eta_0 = 0.001$ were used for AF-STE while $\beta_1 = 0.9$, $\beta_2 = 0.999$, and $\eta_0 = 0.0001$ were used for GEGD. For GEGD, a step size scheduling that increases with the distance from the origin was used from iteration 2 onward to accelerate convergence.
\begin{equation}
\label{eq:step_size}
    \eta^{\left( i \right)} = \eta_0 \\
    \left[ \frac{\left\| \vec{\mu}_\text{L}^{\left( i - 1 \right)} \right\|}{\left\| \vec{\mu}_\text{L}^{\left( 2 \right)} \right\|} \right]^{1/3}
\end{equation}
where the superscript $\left( i \right)$ denotes the iteration number. The latent densities (AF-STE) and the mean latent densities (GEGD) are bounded between $-1$ and 1. Since ADAM does not allow for straightforward inclusion of bounds, we implemented those implicitly using dummy variables $\vec{\zeta}$:
\begin{equation}
\label{eq:dummy_variables}
    \vec{\rho}_\text{L} \\
    = -1 + 2 / \left[ 1 + \text{exp} \left( -\vec{\zeta} \right) \right]
\end{equation}
$\vec{\zeta}$ is unbounded, while $\vec{\rho}_\text{L}$ is asymptotically bound between $-1$ and 1.

A covariance scaling factor of $\sigma_\text{R} = 0.005$ was used for GEGD. In addition, the cost functions for GEGD were exponentiated as $f_\text{exp} = -\text{exp} \left( -\beta_\text{exp} f \right)$, where $\beta_\text{exp} > 1$ controls the exponentiation strength. This exponentiation greatly accentuates the difference between minima with similar depths, allowing the optimizer to more easily overcome its “attraction” towards subpar minima. The transform also aids stable convergence by making gradient directions point more consistently towards promising regions in the parameter space. $\beta_\text{exp} = 20$ was used for all optimization benchmark functions in this work.

\subsection{\label{sec:covariance_regularization}RBF Covariance Regularization}

The regularization factor used in Equation (\ref{eq:RBF_covariance}) is given by:
\begin{equation}
    \epsilon = \frac{\lambda_\text{max} - \kappa \lambda_\text{min}}{\kappa - 1}
\end{equation}
where $\lambda$ are the eigenvalues of the non-regularized $\boldsymbol{\Sigma}_\text{RBF}$ and $\kappa$ is the target condition number after regularization.

\subsection{\label{sec:adaptive_r_CV}Adaptive Determination of the Number of Control Variate Evaluations}

The optimal average variance reduction by approximate control variates depends on the average correlation between the vector functions $\vec{q}f$ and $\vec{q}h$ (Eq.~(\ref{subeq:varianc_reduction_approx})). The correlation is defined under a Gaussian distribution (i.e., different from the global correlation), which means that it may change as the optimization progresses. Therefore, the correlation computed at each iteration determines the optimal $r_\text{CV}$ that produces the optimal variance reduction. The average variance of $\vec{g}_\text{ACV}$ computed with $M$ high-fidelity cost function samples relative to that of $\vec{q}f$ computed with single samples is given by:
\begin{equation}
\label{eq:variance_reduction_adaptive_rCV}
    \frac{\overline{\text{V}}_M \left[ \vec{g}_\text{ACV} \right]}{\overline{\text{V}}_1 \left[ \vec{q}f \right]} \\
    = \frac{1}{M} \left( 1 - \frac{r_\text{CV} - 1}{r_\text{CV}} C^2 \right)
\end{equation}
where $C = \overline{\text{Corr}} \left[ \vec{q}f, \vec{q}h \right]$. $r_\text{CV}$ is constrained by the maximum computation cost per iteration ($t_\text{iter}$).
\begin{equation}
\label{eq:t_iter_inequality}
    M \left( t_\text{HF} + r_\text{CV} t_\text{LF} \right) \leq t_\text{iter}
\end{equation}
Substituting $M = t_\text{iter} / \left( t_\text{HF} + r_\text{CV} t_\text{LF} \right)$ into Eq.~(\ref{eq:variance_reduction_adaptive_rCV}) and solving for the minimizing $r_\text{CV}$ yields:
\begin{equation}
\label{eq:r_CV_float}
    r_\text{CV} = C \sqrt{ \frac{t_\text{HF}}{ t_\text{LF} \left( 1 - C^2 \right)}}
\end{equation}
Note that the $r_\text{CV}$ in Eq.~(\ref{eq:r_CV_float}) is not necessarily an integer. Using this and Eq.~(\ref{eq:t_iter_inequality}), the optimal number of high-fidelity evaluations ($M^*$) and low-to-high evaluation ratio ($r_\text{CV}^*$) are given by the following equations:
\begin{equation}
    M^* = \text{floor} \left[ \frac{t_\text{iter}}{t_\text{HF} + C t_\text{LF} \sqrt{\frac{t_\text{HF}}{t_\text{LF} \left( 1 - C^2 \right)}}} \right]
\end{equation}
\begin{equation}
    r_\text{CV}^* = \text{floor} \left[ \frac{t_\text{iter} - M^* t_\text{HF}}{M^* t_\text{LF}} \right]
\end{equation}
Because the correlation for any given iteration can only be computed after $M^*$ and $r_\text{CV}^*$ are determined, it is instead estimated by the correlation for the previous iteration.

\section{Data Availability}

The code implementation supporting the findings of this paper is available at \url{https://github.com/apmd-lab/gaussian_ensemble_gradient_descent}. The data presented in the figures are available from the corresponding author upon reasonable request.

\section{Conflict of Interest}

The authors declare no conflicts of interest.

\section{Author Contributions}

S.M. conceived the original idea and developed the theory with input from J.P. and J.S.; S.M. implemented the algorithms and ran benchmark tests; J.S. supervised the project. The manuscript was written by S.M. with input from J.S.

\begin{acknowledgments}
This work is supported by the National Research Foundation of the Republic of Korea (RS-2024-00414119, RS-2023-00283667).
\end{acknowledgments}

\bibliography{references}
\end{document}



\title{Supplementary Information for\\Ensemble-Based Global Search Framework for the Design Optimization of Fabrication-Constrained Freeform Devices}


\author{Seokhwan Min}
\author{Junhyung Park}
\author{Jonghwa Shin}
\affiliation{%
 Department of Materials Science and Engineering,\\
 Korea Advanced Institute of Science and Technology,\\
 Daejeon 34141, Republic of Korea
}%


\maketitle

\newpage

\section{\label{sec:1}S1. Test Function Details}

\setcounter{figure}{0}
\renewcommand\thefigure{S\arabic{figure}}

\begin{figure}[h]
\includegraphics[width=0.9\textwidth]{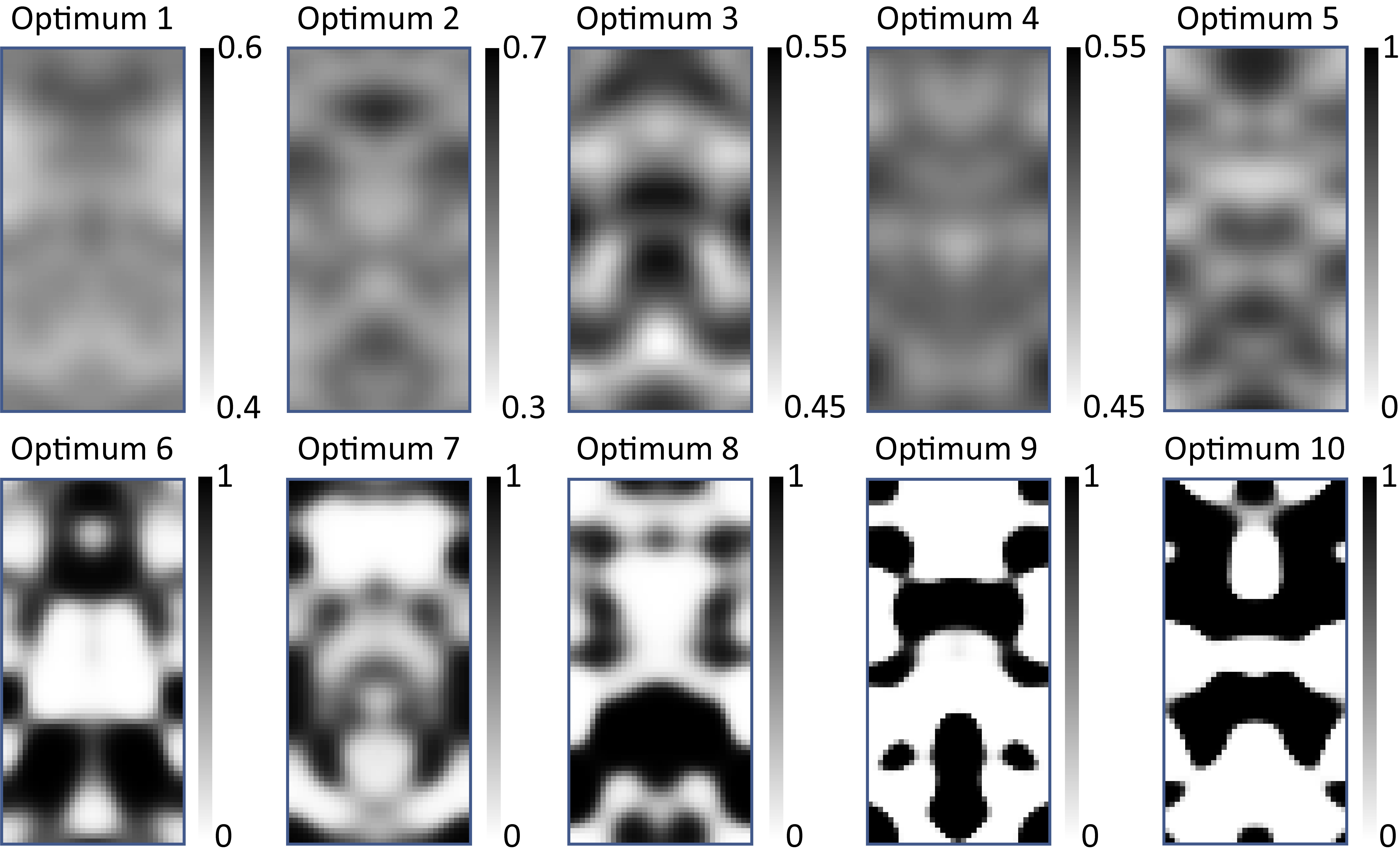}
\caption{\label{fig:test_function} The locally optimal grayscale designs used to define the test function in the main text.}
\end{figure}

\begin{figure}[h]
\includegraphics[width=0.4\textwidth]{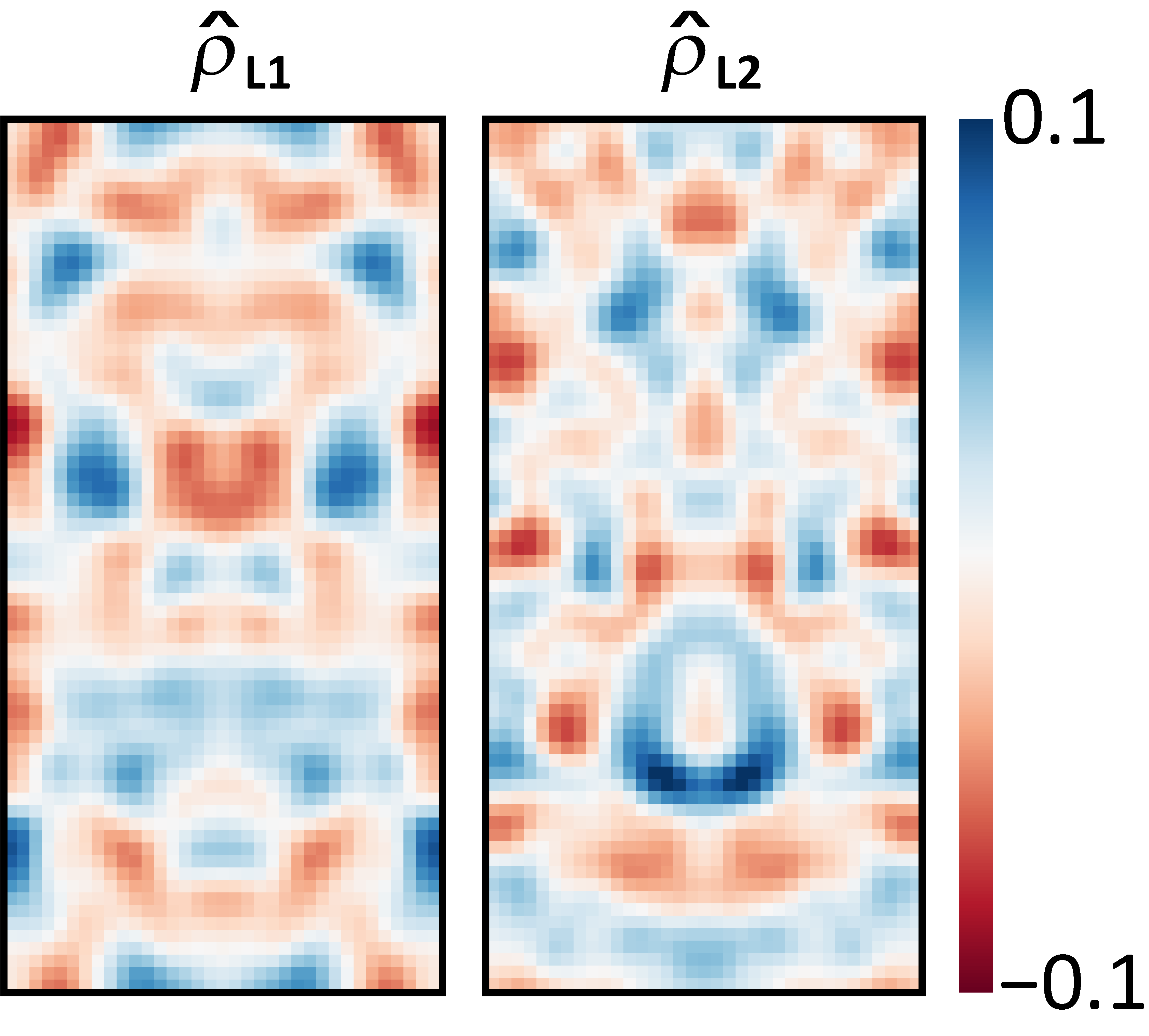}
\caption{\label{fig:2d_cross_section_unit_vectors} 2D representation of the unit vectors along which the cross section of the test function was computed in Figure 5a of the main text.}
\end{figure}

\clearpage

\section{\label{sec:2}S2. Benchmark Tests on the Effectiveness of the Sampling Efficiency Enhancement Techniques\protect}

Here we benchmark the effectiveness of the various sampling efficiency enhancement strategies for the test function used in the main text. The following CV was defined as a low-fidelity model for the test function:
\begin{equation}
    f_\text{test,CV} = f_\text{test} + 0.001 \eta \left( \vec{\rho}_\text{F} \right)
\end{equation}
where $\eta \sim N \left( 0, 1 \right)$ is a fixed Gaussian noise dependent on $\vec{\rho}_\text{F}$ (i.e., different $\vec{\rho}_\text{F}$ will have different $\eta$, but $\eta$ is fixed throughout optimization). A total of $M = 10$ high-fidelity samples were used for optimization runs without CV, and a total of $M \geq 5$ high-fidelity samples with an adaptive control variate sampling ratio $r_\text{CV}$ was used for runs with CV. $M$ and $r_\text{CV}$ were determined adaptively each iteration based on the correlation between the low and high-fidelity samples from the previous iteration. For every iteration, it was ensured that $M \times t_\text{HF} + M \times r_\text{CV} \times t_\text{LF} \approx 10 t_\text{HF}$ where $t_\text{HF}$ and $t_\text{LF}$ are runtimes for the high and low-fidelity simulations, respectively.

\begin{figure}
\includegraphics[width=0.5\textwidth]{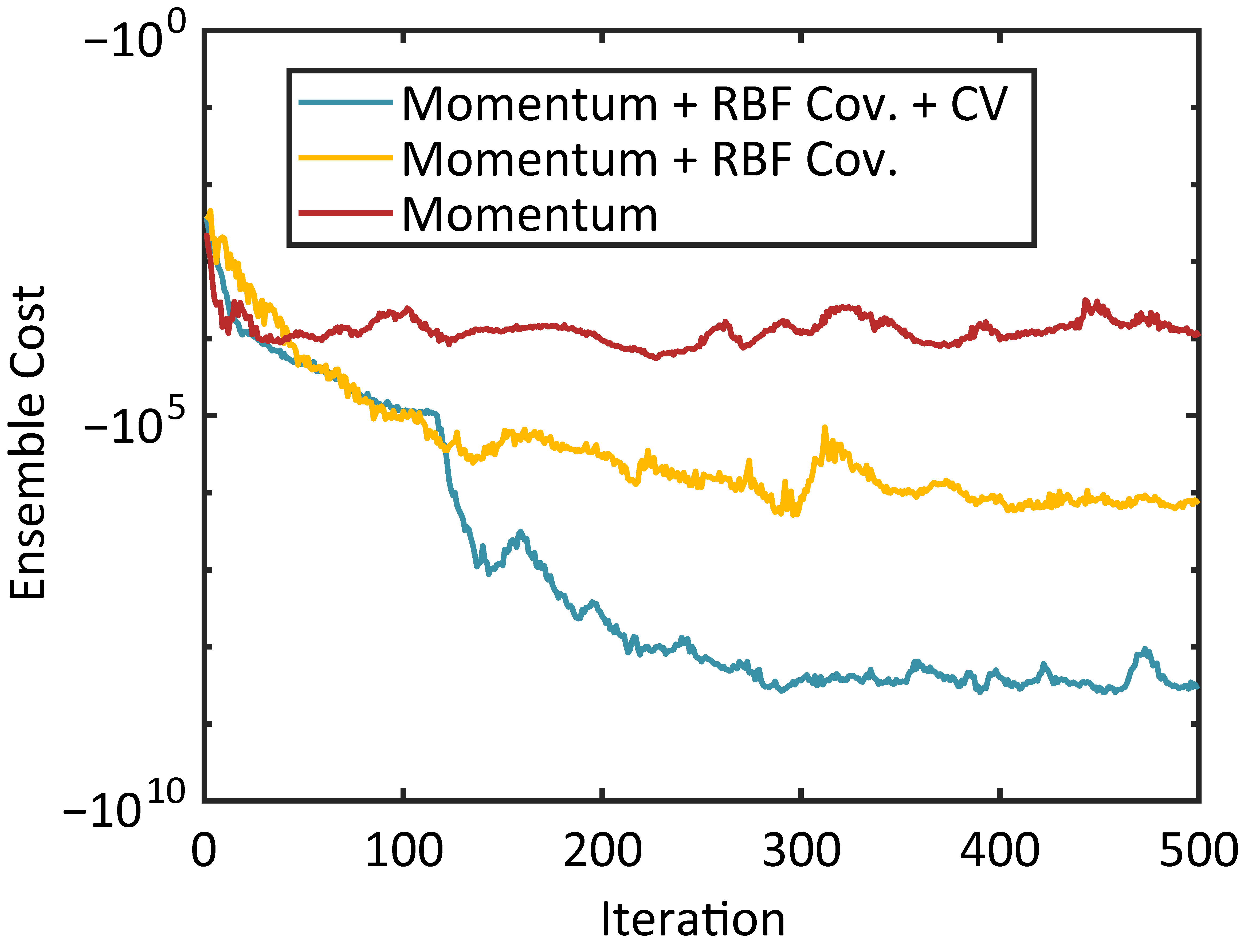}
\caption{\label{fig:sampling_efficiency_enhancements} Ensemble cost evolution during GEGD using different levels of sampling efficiency enhancements.}
\end{figure}

The effectiveness of the sampling efficiency enhancement strategies was gauged by comparing the ensemble average cost evolution instead of the best cost obtained, as the best cost depends heavily on random chance and does not accurately reflect the actual quality of optimization. As can be seen in Figure~\ref{fig:sampling_efficiency_enhancements}, the application of each method (RBF covariance and CV) significantly improves the converged ensemble average cost, justifying their inclusion in GEGD despite the added complexity.

\clearpage

\section{\label{sec:3}S3. Optimization Trajectories and Optimized Designs for the Nanophotonic Benchmark Problems\protect}

\begin{figure}[h]
\includegraphics[width=\textwidth]{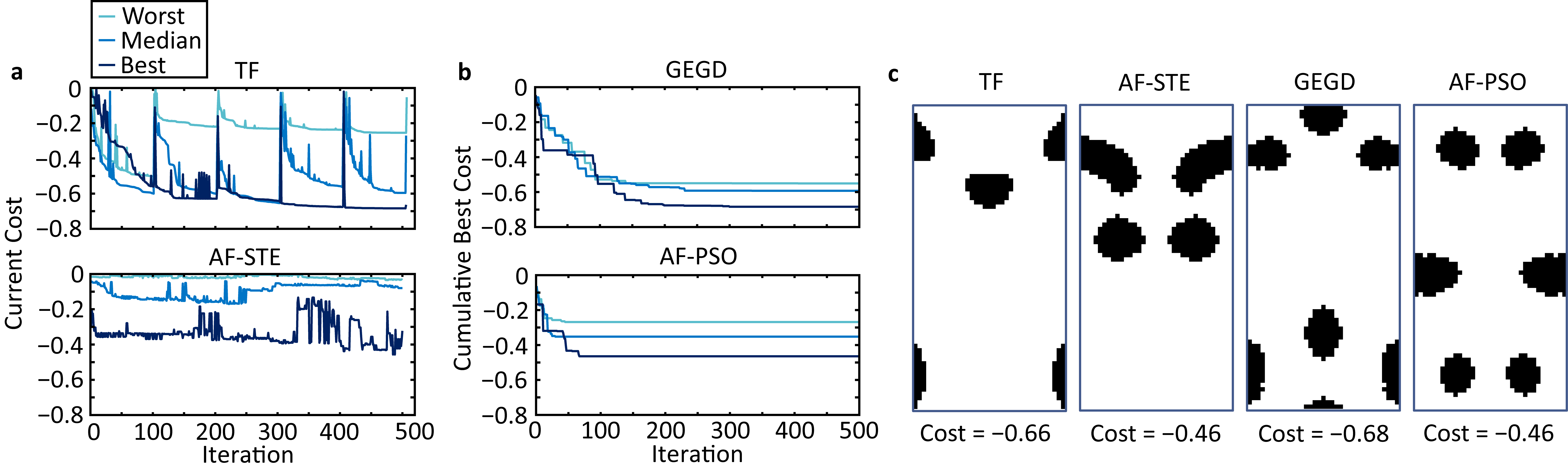}
\caption{\label{fig:polarization_beamsplitter_results} Optimization benchmark results for the polarization beamsplitter. \textbf{a}, Cost evolution during optimization for selected TF and AF-STE runs with the worst, median, and best final costs. \textbf{b}, Cumulative best cost evolution during optimization for selected GEGD and AF-PSO runs with the worst, median, and best final costs. \textbf{c}, Best final designs for each benchmark algorithm.}
\end{figure}

\begin{figure}[h]
\includegraphics[width=\textwidth]{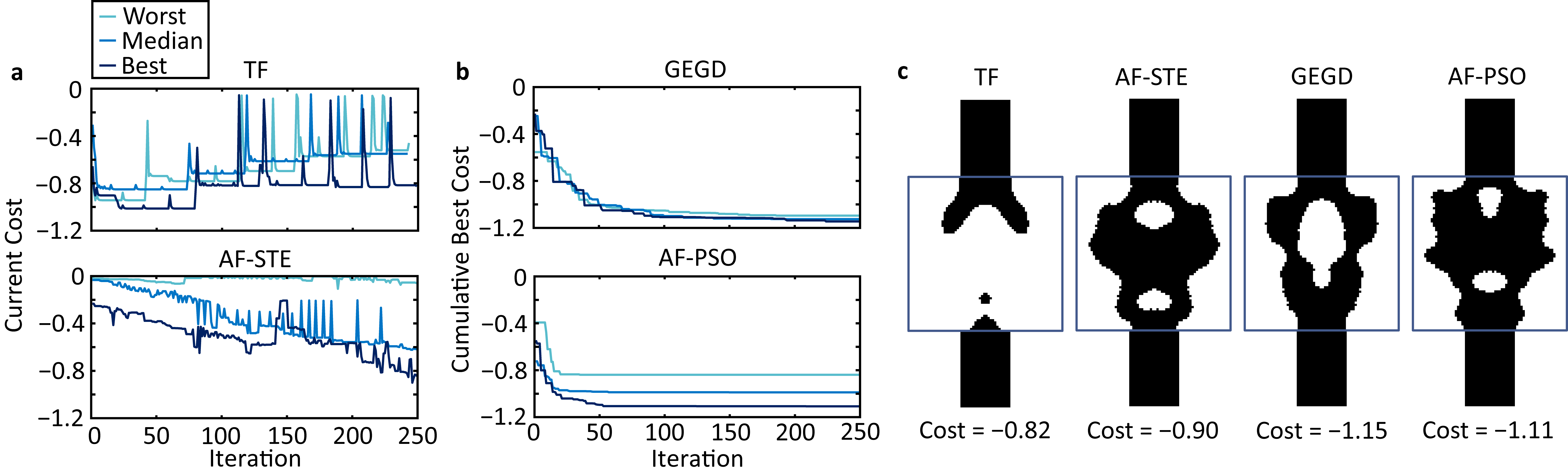}
\caption{\label{fig:mode_converter_results} Optimization benchmark results for the 2D mode converter. \textbf{a}, Cost evolution during optimization for selected TF and AF-STE runs with the worst, median, and best final costs. \textbf{b}, Cumulative best cost evolution during optimization for selected GEGD and AF-PSO runs with the worst, median, and best final costs. \textbf{c}, Best final designs for each benchmark algorithm.}
\end{figure}
